\documentclass[10pt]{article}

\usepackage{amsmath}
\usepackage{amssymb}

\usepackage[all]{xy} 

\usepackage{pstricks}

\newcounter{remark}
\newcommand{\remark}{\addtocounter{remark}{1}
                       \par \quad {\bf \arabic{remark}}.\,
                      }

\newenvironment{proof}{\medbreak \noindent {\bf Proof~---}}
                       {\hfill $\square$ \medbreak}

\newcommand{\Cl}{\operatorname{Cl}}

\newcommand{\Gal}{\operatorname{Gal}}

\newcommand{\ZZ}{\mathbb Z}
\newcommand{\QQ}{\mathbb Q}
\newcommand{\RR}{\mathbb R}
\newcommand{\CC}{\mathbb C}
\newcommand{\FF}{\mathbb F}

\newcommand{\PP}{\mathbb P}

\newcommand{\OO}{\mathcal O}
\newcommand{\HH}{\mathcal H}
\newcommand{\XX}{\mathcal X}

\newtheorem{theo}{Theorem}

\newtheorem{cor}[theo]{Corollary}

\newtheorem{prop}[theo]{Proposition}

\title{\bfseries Computation of the cover of Shimura curves~$\XX_0(2) \to \XX(1)$ for the cyclic cubic field
of discriminant~$13^2$}

\author{
Emmanuel
{\sc Hallouin}\thanks{Institut de Mathématiques de Toulouse (Laboratoire \'Emile Picard),
University of Toulouse~II, France.}
}

\begin{document}
\maketitle

\begin{quote}
  \footnotesize {\bfseries Summary ---} We compute the canonical model of the cover of Shimura
curves $\XX_0(2) \to \XX(1)$ for the cubic field of discriminant~$13^2$ described at the end
of~Elkies' paper~\cite{Elkies}. Last, we list the coordinates of some rational CM points
on~$\XX(1)$.
\end{quote}


\section*{Introduction}

In~\cite{Elkies}, Noam D. Elkies computes explicit algebraic models of some Shimura
curves associated with quaternion algebras whose center is a cubic cyclic totally real field
and whose ramification data consists in exactly two of the infinite places of the center.
The examples of~$K_7$ and~$K_9$, the cubic fields of discriminant~$49 = 7^2$ and~$81 = 9^2$
respectively, are studied in details. Then the case of the center~$K = K_{13}$, the cubic field
of discriminant~$13^2$, is described. Elkies shows that~$\XX(1)$ has genus zero and that the
cover~$\XX_0(2) \to \XX(1)$ is a degree~$9$ cover with four branch points
and with Galois group equal to~$PSL_2(\FF_8)$. But the explicit computation of the canonical model
remains to be done.

It turns out that this cover belongs to a family of degree~$9$ covers of~$\PP^1_\CC$ with
Galois group~$PSL_2(\FF_8)$, that I have computed in~\cite{PSL_2-F_8}. The strategy followed
to compute the cover~$\XX_0(2) \to \XX(1)$ is though very simple. First gathering enough
data's
about the cover of Shimura curves (see~\S\ref{s_The_Curve}), secondly putting these conditions
into equations in order to compute the parameter corresponding to the expected element
(see~\S\ref{s_Family}).

One of the key properties of~$\XX_0(2)$ is that it has a an involution~$w_2$.
Once we know explicitly the map~$\XX_0(2) \to \XX(1)$ and this involution~$w_2$,
it is an easy task to deduce the coordinates of some CM points on~$\XX(1)$. We obtain,
this way, some rational CM points on~$\XX(1)$.

\section{The cover of Shimura curves~$\XX_0(2) \to \XX(1)$} \label{s_The_Curve}

Let~$K$ be the unique cyclic cubic field of discriminant~$13^2$, which can be obtained by adjoining
to~$\QQ$ a root~$\theta_0$ of~$x^3 + x^2 - 4x + 1 \in \QQ[x]$. Its narrow class group~$\Cl^+(K)$,
and thus its class group~$\Cl(K)$, are trivial. We denote
by~$\sigma_i$,~$1 \leq i \leq 3$, the three real embeddings of~$K$ into~$\RR$.

Let~$B$ be the unique (up to $K$-isomorphisms) quaternion algebra over~$K$ ramified
at~$\sigma_2$ and~$\sigma_3$ and nowhere else.  The embedding~$\sigma_1$ being unramified in~$B$,
we know that~$B \otimes_{K, \sigma_1} \RR \simeq M_2(\RR)$ (the notation~${} \otimes_{K, \sigma_1} \RR$ means that~$\RR$ is considered as a $K$-algebra via the
embedding~$K \overset{\sigma_1}{\hookrightarrow} \RR$). This allows us to fix an
embedding~$i_{\infty} : B \hookrightarrow M_2(\RR)$ which restricts to~$\sigma_1$ on~$K$.
We also fix~$\OO$ one of the maximal orders of~$B$;
since~$\Cl^+(K)$ is trivial, all these maximal orders are known to be conjugate. We
denote by~$\OO^1$ the set of invertible elements of~$\OO$ whose reduce norm equals~$1$.
Then the group~$\Gamma(1) \overset{\text{def.}}{=} i_\infty(\OO^1)$ is an arithmetic Fuchsian
group (cf.~\cite{Katok}), and in particular a discrete subgroup of~$PSL_2(\RR)$ (I denote the same
way a subgroup of~$SL_2(\RR)$ and its image in the quotient~$PSL_2(\RR)$).

Let~$\HH$ be the usual upper half-plane with the action of~$PSL_2(\RR)$.
The quotient~$\Gamma(1) \setminus \HH$, denoted by~$\XX(1)_\CC$, is a compact Riemann surface, called
the Shimura curve associated to~$B$. It has a structure of complex algebraic curve.

The (inert) prime~$2$ in~$K$ being unramified in~$B$, we know
that~$B \otimes_K \QQ_8 \simeq M_2(\QQ_8)$. We can fix an
embedding~$i_2 : B \hookrightarrow M_2(\QQ_8)$ in such a way that~$i_2(\OO) = M_2(\ZZ_8)$.
As for the classical modular curves case, we consider the subgroups~$\Gamma_0(2), \Gamma_1(2)$
and~$\Gamma(2)$ of~$\Gamma(1)$. For example,~$\Gamma_0(2)$, which is the subgroup we will focus on,
 is the image by~$i_\infty$
of:
$$
\left\{\omega \in \OO^1, \, i_2(\omega) \equiv \begin{pmatrix}* & *\\0 & *\end{pmatrix}
\bmod{2M_2(\ZZ_8)}\right\}.
$$
They are also discrete subgroups of~$PSL_2(\RR)$ and the quotients
of~$\HH$ by these subgroups, respectively denoted by~$\XX_0(2)_\CC, \, \XX_1(2)_\CC$
and~$\XX(2)_\CC$, are complex algebraic curves. We have the tower:
$$
\xymatrix{
\XX(2)_\CC \ar@{->}[r]^8
\ar@/^2pc/@{-}[rrr]|{PSL_2(\FF_8)} &
\XX_1(2)_\CC \ar@{->}[r]^7 &
\XX_0(2)_\CC \ar@{->}[r]^9 &
\XX(1)_\CC.
}
$$
The cover~$\XX(2)_\CC \to \XX(1)_\CC$ is the Galois closure of the
cover~$\XX_0(2)_\CC \to \XX(1)_\CC$.
Moreover, the curve~$\XX_0(2)_\CC$ has an involution~$w_2$ which comes from a trace-zero element
of~$B$ of norm~$2$.

More precisely, in~\cite{Elkies} pages~$314$-$315$, thanks to Shimizu's area formula, Elkies proves
that~$\Gamma(1)$ has four elliptic points, one of
order~$3$,~$P_0$, three of order~$2$,~$P_1$, $P_2$, $P_3$ and that the curve~$\XX(1)$ has genus zero.
The map~$\XX_0(2)_\CC \to \XX(1)_\CC$ is a cover of degree~$9 = \#\FF_8 + 1$, with geometric Galois
group~$PSL_2(\FF_8)$ (acting on the nine points of~$\PP^1_{\FF_8}$) and is ramified at the elliptic
points~$P_0,P_1,P_2,P_3$. Since elements of order~$3$ and~$2$ in~$PSL_2(\FF_8)$ have cycle
shapes~$3^3$ and~$2^3 \cdot 1$ respectively, the ramification data of this cover must
be~$(3^3, 2^4 \cdot 1, 2^4 \cdot 1, 2^4 \cdot 1)$.
We denote by~$Q_1, Q_2, Q_3$ the there unramified
points on~$\XX_0(2)_\CC$ above~$P_1,P_2,P_3$ respectively.

Apart form elliptic points, Shimura curves also contains CM-points who are associated to CM-fields
over~$K$.
Let~$L/K$ be a CM-field and~$\OO_L$ its ring of integers. The field~$L$ can be embed
in~$B$ by~$i : L \hookrightarrow B$. The ring~$i^{-1}(i(L) \cap \OO)$ is an order of~$\OO_L$
and the embedding~$i$ is called optimal if and only this order is~$\OO_L$ itself.
Then all the elements of~$i_\infty \circ i(L)$
share the same fixed point~$\tau \in \HH$. The point~$\tau \bmod{\Gamma(1)} \in \XX(1)_\CC$ is called
a CM-point on~$\XX(1)_\CC$ by the CM-order~$i^{-1}(i(L) \cap \OO)$ or by~$L$ if the embedding~$i$
is optimal.
 One can show that two
embeddings~$i$ and~$i'$, whose associated orders are equal, give rise to the same CM-point
on~$\XX(1)_\CC$ if and only if they
are conjugate by a invertible element of~$\OO$. Since there are finitely many embeddings
up to conjugacy, there are finitely many CM-points on~$\XX(1)_\CC$ by a given order.
In fact their number is just the class number of the order.

\medbreak

{\bfseries\noindent Shimura's results.}
Thanks to a modular interpretation, Shimura proved that the curve~$\XX(1)_\CC$ have a canonical
model defined over the narrow class field of~$K$, namely~$K$ itself here.

\begin{theo}[Shimura,~\cite{Shimura}, \S3.2, Main Theorem~I]
There exists an algebraic curve~$\XX(1)$ and a holomorphic
map~$j : \HH \to \XX(1)$ satisfying the following conditions.
\begin{enumerate}
\item The curve~$\XX(1)$ is defined over~$K$.
\item The map~$j$ yields an analytic isomorphism from~$\Gamma(1) \setminus \HH$
to~$\XX(1) \otimes_K \CC$.
\item Let~$L$ be a CM-field over~$K$ which optimally embeds in~$B$ and whose class number is denoted
by~$h_L$. Then there are exactly $h_L$ CM-points~$\tau_1, \ldots, \tau_{h_L}$ by~$L$
on~$\XX(1)$. The values~$j(\tau_i)$, for~$1 \leq i \leq h_L$, form a complete set of conjugates
over~$K$ and for each~$i$, one has~$L^{\text{hilb.}} = L \cdot K(j(\tau_i))$.
\end{enumerate}
\end{theo}

As in the classical case, computing the modular polynomial relying the functions~$j$
and~$j \circ w_2$, one can prove that the
curve~$\XX_0(2)$ is also defined over~$K$.

As for the curve~$\XX(2)$, Shimura proves that it
admits a canonical model which is defined over the strict $2$-ray class field of~$K$
(i.e. unramified outside~$2$ and
the real places), namely~$K$ another time since~$\Cl^+(K, 2)$ is also trivial.

\medbreak

{\bfseries\noindent Galois descent to~$\QQ$.}
In fact, as in~\cite{Elkies_SCC} (page~38) or~\cite{Voight}\S6, Galois descent to~$\QQ$ is
possible.
Nevertheless, the proofs of this fact contained in these two papers should be adapted
since, here, the group~$\Gamma(1)$ is no more a triangle group.

The $K$-model of~$\XX(1)_\CC$ is strongly connected to the CM-points; the elliptic points play
a crucial role in the fact that this model also descents to~$\QQ$.

The elliptic point~$P_0$ of order~$3$ is CM
by~$K(\sqrt{-3})$ which has class number~$1$. Thanks to Shimura's result,
the point~$P_0$ is defined over~$K$.

The elliptic points~$P_1, P_2, P_3$ of order~$2$ are CM by~$K(\sqrt{-1})$ which has class number
equal to~$3$. Due to Shimura's result, these points must be conjugate to  each other over~$K$.


\begin{prop}
The curve~$\XX(1)$ is defined over~$\QQ$, in such a way that:
\begin{enumerate}
\item the unique elliptic point~$P_0$ of order~$3$, which is
CM by~$K(\sqrt{-3})$, is rational over~$\QQ$;
\item the three elliptic points~$P_1,P_2,P_3$ of order~$2$, who are CM by~$K(\sqrt{-1})$,
are conjugate to each others over~$\QQ$.
\end{enumerate}
\end{prop}

\begin{proof}
The proposition follows from the fact that, not only the complete curve~$\XX(1)$,
but also the curve~$\XX(1) \setminus \{P_0\} \cup \{P_1,P_2,P_3\}$ descents to~$\QQ$
(the point~$P_0$ is ``colored'', while the three others~$P_1,P_2,P_3$ are ``uncolored'').
The key point is that the cubic field~$K$, we start with, is a Galois extension of~$\QQ$.

By functoriality, for each~$\sigma \in \Gal(K/\QQ)$, the curve~$\XX(1)^\sigma$ is nothing else that the
Shimura curve corresponding to the algebra~$B^\sigma$ which is ramified at~$\sigma_2 \circ \sigma$
and~$\sigma_3 \circ \sigma$ and nowhere else.

The algebras~$B$ and~$B^\sigma$ are $\sigma$-isomorphic, meaning that there exists
a~$\QQ$-isomorphism~$\theta_\sigma : B \to B^\sigma$ which restricts to~$\sigma$ on~$K$.
Let us denote by~$\Gamma(1)_\sigma$ (I put~$\sigma$ in subscript to emphasize the
fact~$\Gamma(1)_\sigma$ is
not the group~$\Gamma(1)$ whose coefficient are twisted by~$\sigma$; this would not make any sense)
the arithmetic group corresponding to~$B^\sigma$; it
follows form the choice of a maximal order~$\OO^\sigma$ of~$B^\sigma$ and an embedding
of~$B^\sigma$ in~$M_2(\RR)$. Since~$\theta_\sigma^{-1}(\OO^\sigma)$ is also a maximal order of~$B$,
it must be conjugate to~$\OO$ by an element of~$B$. Concerning the groups~$\Gamma(1)$
and~$\Gamma(1)_\sigma$,
this means that there exists~$\gamma \in PSL_2(\RR)$ such
that~$\Gamma(1)_\sigma = \gamma \Gamma(1) \gamma^{-1}$. Then there exists an
isomorphism of Riemann surfaces between~$\XX(1)_\CC$ and~$\XX(1)^\sigma_\CC$ which yields
an isomorphism~$\varphi_\sigma : \XX(1) \to \XX(1)^\sigma$ of algebraic curves over~$K$.

It is easy to verify that these isomorphisms map an elliptic point of a given order to an
elliptic point of the same order. So they are isomorphisms between the
curve~$\XX(1) \setminus \{P_0\} \cup \{P_1, P_2, P_3\}$ and its
twist~$\XX(1)^\sigma \setminus \{{}^\sigma P_0\} \cup \{{}^\sigma P_1, {}^\sigma P_2, {}^\sigma P_3\}$, for~$\sigma \in \Gal(K/\QQ)$.

To show that the~$\varphi_\sigma$'s give a sufficient data for the Galois descent to~$\QQ$, it
suffices to prove that the curve~$\XX(1) \setminus \{P_0\} \cup \{P_1,P_2,P_3\}$ and its
conjugates do not
have any non trivial automorphism. Such an automorphism is an automorphism of genus zero curve
which must fix the point~$P_0$ and
permute~$P_1,P_2,P_3$. So its order is~$2$ or~$3$, depending on the fact that it fixes or not one
of the~$P_i$'s for~$i \geq 1$. The cross ratios~$[P_0,P_1,P_2,P_3]$ is then equal
to~$[\infty, 0, 1, -1]$ or~$[\infty, 1, j, j^2]$ (where~$j = e^{\frac{2i\pi}{3}}$). On the
other hand, since~$P_0$ is rational over~$K$ and since the~$P_i$'s for~$i \geq 1$ are conjugate to
each others over~$K$ and such that~$K(\sqrt{-1}) \cdot K(P_i) = K(\sqrt{-1})^{\text{hilb.}}$, the cross
ratio~$[P_0,P_1,P_2,P_3]$
must satisfy:
$$
K([P_0,P_1,P_2,P_3]) \cdot K(\sqrt{-1}) = K(\sqrt{-1})^{\text{hilb.}}.
$$
Neither~$[\infty, 0, 1, -1]$ nor~$[\infty, 1, j, j^2]$ satisfies this, therefore the
curve~$\XX(1) \setminus \{P_0\} \cup \{P_1,P_2,P_3\}$ do not have any non trivial
automorphism.

In this case, the Weil descent criterion is automatically satisfied.
\end{proof}


\begin{prop}
The curve~$\XX_0(2)$ is defined over~$\QQ$ and the three elliptic points $Q_1$, $Q_2$, $Q_3$ of
order~$2$ are conjugate to each others over~$\QQ$ and fixed by the involution~$w_2$.
\end{prop}

\begin{proof}
The key point to prove the decent is that the prime~$2$ of~$K$ is $\Gal(K/\QQ)$-invariant.
Indeed, keeping the notation of the preceding proof, this permit to show that the
groups~$\Gamma_0(2)$ and~$\Gamma_0(2)_\sigma$ are also
conjugate to each other by the same~$\gamma$. Then the isomorphism~$\varphi_\sigma$ from~$\XX(1)$
to~$\XX(1)^\sigma$ lifts to an isomorphism form~$\XX_0(2)$ to~$\XX_0(2)^\sigma$.

On the other hand, since~$PSL_2(\FF_8)$ is self-centralizing in the group of permutations of the
points of~$\PP^1_{\FF_8}$, the cover~$\XX_0(2) \to \XX(1)$ and its conjugates do not have any non
trivial automorphism. Therefore the previous isomorphisms between~$\XX_0(2)$ and its conjugate
give a sufficient conditions to the descent to~$\QQ$.

The construction of the involution~$w_2$ being Galois-invariant, this function must be defined
over~$\QQ$. Moreover, since~$w_2$ permutes elliptic points of same order, it permutes
the~$Q_i$'s. So~$w_2$ must fix one of them. By Galois conjugation, it fixes all of them.
\end{proof}

\section{Selecting the good element in the family} \label{s_Family}

In~\cite{PSL_2-F_8}, I have computed an algebraic model of
the universal family of degree~$9$ covers of~$\PP^1_\CC$ with Galois group~$PSL_2(\FF_8)$,
with ramification type~$(3^3, 2^4 \cdot 1, 2^4 \cdot 1, 2^4 \cdot 1)$ and such that the first
branch point is rational while the three others are conjugate to each other.
The final result (proposition~9) consists in a total degree~$3$ polynomial of~$\QQ(T)[f,g]$
which is an equation of a genus~$1$ curve~$\mathcal{E}$ over~$\QQ(T)$,
and a function~$\varphi \in \QQ(T)(\mathcal{E})$. This function is such that the
map~$\mathcal{E} \overset{\varphi}{\to} \PP^1_{\QQ(T)}$ is a degree~$9$ cover,
with (geometric) Galois group~$PSL_2(\FF_8)$, with four branch points: the point of
type~$3^3$ is at~$\varphi = \infty$, the three points of type~$2^4 \cdot 1$ have
$\varphi$-coordinate equal to the three roots
of~$x^3 + H(T)(x+1)$, where~$H \in \QQ(T)$ is explicit. The divisor of the degree~$3$ function~$f$
is related to the ramification points as follows: its zeros are three unramified points
on~$\mathcal{E}$ over the three branch points of type~$2^4 \cdot 1$ while its poles are the
three ramification points of index~$3$.

The variable~$T$ is a coordinate on~$\HH$, the (absolute) Hurwitz space parameterizing this family,
which thus is~$\QQ$-isomorphic to~$\PP^1_\QQ$. Besides this space, I also have computed the (inner)
Hurwitz space~$\HH_{PSL_2(\FF_8)}$ parameterizing the $PSL_2(\FF_8)$-Galois covers. It is
$\QQ$-isomorphic to~$\PP^1_\QQ$, with parameter~$S$, and it covers~$\HH$ by an explicit degree~$3$
map (see proposition~10, loc. cit.).

The specializations of this family at rational values of the parameter~$T$ yield $PSL_2(\FF_8)$
covers of~$\PP^1_\QQ$ defined over~$\QQ$, with four branch points and expected ramification data.

It turns out that the cover~$\XX_0(2) \to \XX(1)$ is such a specialization.
Which one? That is the question! 

\begin{prop}
The cover~$\XX_0(2) \to \XX(1)$ corresponds to the specialization at the value~$T = -1$ of
our family.
\end{prop}

\begin{proof}
First of all, the cover~$\XX_0(2) \to \XX(1)$ is truly a member of our family. Indeed the
ramification data, the Galois group do correspond and we know that the unique elliptic point
of order three on~$\XX(1)$ is rational and that the three others of order~$2$ are conjugate
to each others.

Let~$\widetilde{\mathcal{E}}$ be the specialization we are looking for. We know that~$\widetilde{\mathcal{E}}$ is
defined over~$\QQ$ and that there exists an involution~$w_2$ of~$\widetilde{\mathcal{E}}$ also defined
over~$\QQ$ that fix the three unramified points~$Q_1,Q_2,Q_3$ over~$P_1,P_2,P_3$
the three elliptic points of~$\XX(1)$ of order~$2$. Extending the scalars to~$\QQ(Q_1)$, the
genus~$1$ curve~$\widetilde{\mathcal{E}}$ becomes an elliptic curve whose origin can be chosen to
be equal to~$Q_1$. Then~$w_2$ is nothing else that~$P \mapsto -P$ and the $2$-torsion consists in
the~$Q_i$'s plus the point~$Q_2 + Q_3$. Therefore one has the relations~$2Q_1 = 2Q_2 = 2Q_3$.

These conditions completely determine the parameter~$T$ corresponding to the
cover~$\XX_0(2) \to \XX(1)$. To show this fact, let us go back on the generic curve~$\mathcal{E}$.
We choose the same notation for the
unramified points above the branch points of type~$2^4 \cdot 1$, namely~$Q_1,Q_2,Q_3$.
They are conjugate to each others over~$\QQ(T)$ since they are the zeros of the degree-$3$
function~$f \in \QQ(T)(\mathcal{E})$. We extend the scalars to~$\QQ(T)(\Theta)$ the splitting field
of~$\QQ(T)(Q_i)$ the field of definition of the~$Q_i$'s; this is a degree~$6$ extension
over~$\QQ(T)$. We choose the
point~$Q_1$ as the origin of~$\mathcal{E}$ over~$\QQ(T)(\Theta)$ and we compute the
$(f,g)$-coordinates of the points~$2Q_2$ and~$2Q_3$; they are unfortunately to huge to figure in
this paper.
The parameter~$T_0 \in \QQ$, we are looking for, must be such that these two points
specialized at~$T_0$ are equal. This leads to many polynomials in~$T$
who must have~$T_0$ as a common zero. But the gcd of these polynomials is equal to~$T+1$.
\end{proof}

Since we know an algebraic model of the universal family~$\varphi : \mathcal{E} \to \PP^1_{\QQ(T)}$
(see~\cite{PSL_2-F_8} \S3.4, proposition~9), it is very easy to compute the
specialization~$\varphi_{-1} : \mathcal{E}_{-1} \to \PP^1_\QQ$.
We deduce an explicit algebraic model of~$\XX_0(2) \to \XX(1)$. Elkies notes that the equation
simplifies somewhat by taking~$(43754\varphi_{-1} + 687)/2061$ in place of~$\varphi_{-1}$.

\begin{cor}
The curve~$\XX_0(2)$ has equation:
$$\textstyle
f^3 - f^2g - f^2 - \frac{25}{84}fg + \frac{40}{147}f - \frac{625}{702}g^3 + \frac{50}{441}g - \frac{640}{9261} = 0
$$
and the cover~$\XX_0(2) \to \XX(1)$ is given by the function:
$$\textstyle
\varphi =
-\frac{241129}{125}f^2g - \frac{36309}{125}f^2 + \frac{1715}{3}fg^2 + \frac{22099}{30}fg - \frac{637}{200}f - \frac{1715}{9}g^3 + \frac{1225}{36}g^2 - \frac{1708}{45}g - \frac{1301}{75}.
$$
It is a degree~$9$ cover ramified at~$\varphi = \infty$ and at the three
roots
of~$x^3 - x^2 - 992x - 20736$.
The  ramification data
is~$(3^3, 2^4 \cdot 1, 2^4 \cdot 1, 2^4 \cdot 1)$ and the (geometric) Galois group
is equal to~$PSL_2(\FF_8)$.
\end{cor}

Knowing this model, one can now check some easy facts.

$\bullet$ Since~$w_2$ is defined over~$\QQ$, the fix points
divisor~$Q_1 + Q_2 + Q_3 + [Q_2+Q_3]$ is rational;  but the~$Q_i$'s being conjugate
to each others, necessarily the point~$Q_2+Q_3$ must be rational. It does and it
has $(f,g)$-coordinates equal to~$\left(\frac{16}{21}, 0\right)$. Thus~$\XX_0(2)$
is an elliptic curve which, as suggested by Watkins (see \S4 in~\cite{Elkies}), is
$\QQ$-isomorphic to the elliptic curve~$[1,-1,1,-65773,-6478507]$.

$\bullet$ The involution~$w_2$ of~$\QQ(f,g)$ is given by:
\begin{align*}
w_2(f) &=
\frac{733824f^2 - 733824fg - 454272f + 1715000g^3 - 218400g + 133120}{4501875g^3 + 733824g + 279552}, \\
w_2(g) &=
\frac{-91728f^2 g + 91728f g^2 + 56784fg + 62244g^2 - 3328g}{214375g^3 + 34944g + 13312}.
\end{align*}
It does fix the~$Q_i$'s and the rational point~$\left(\frac{16}{21}, 0\right)$.

$\bullet$ Last, we recall that the Galois closure~$\XX(2)$ of the cover~$\XX_0(2) \to \XX(1)$ is
known to be defined over~$K$
(since~$\Cl^+(K,2)$ is trivial). Therefore~$K$ must be the field of definition of the pre-image
of the point~$T = -1$ by the cover~$\HH_{PSL_2(\FF_8)} \to \HH$; hopefully, it does!
As in the context of triangle groups where J.Voight has shown in the claim
following theorem~6.2 of~\cite{Voight}, that the curve~$\XX(2)$ may have a model defined
over~$\QQ$. Nevertheless, the automorphisms of the Galois cover~$\XX(2) \to \XX(1)$ are certainly
not defined over~$\QQ$.

\section{Rational CM-points on~$\XX(1)$} \label{s_CM_points}

Before listing rational CM points, we come back over the three elliptic
points~$P_1$, $P_2$, $P_3$ of order~$2$ on~$\XX(1)$. They are CM by~$K(\sqrt{-1})$ and
thanks to Shimura's result, we must have:
$$
K(\sqrt{-1})^{\text{hilb.}} = K(\sqrt{-1}) \cdot K(P_1).
$$
Hopefully, this is true. Indeed, the field of definition of the~$P_i$'s is the cubic
field~$\QQ(\theta_1)$ with~$\theta_1^3 - \theta_1^2 - 4\theta_1 + 12$ whose discriminant is equal
to~$2^2 13^2$ and one can check that~$K(\sqrt{-1})^{\text{hilb.}} = K(\sqrt{-1},\theta_1)$.

Following Elkies, the rational CM points on~$\XX(1)$ must come from a CM-fields~$K(\sqrt{-D})$
for~$D \in \QQ^*_+$ such that~$\QQ(\sqrt{-D})$ is principal. An easy computation shows that
only~$D = 2, 3, 7$ can appear.

\begin{prop}
On~$\XX(1)$ there are
CM points who are defined over~$\QQ$:
\begin{enumerate}
\item the pole of~$\varphi$, namely the unique elliptic point of order~$3$, which is CM
by~$K(\sqrt{-3})$,
\item the point with $\varphi$-coordinate equal to~$\frac{2^4 \cdot 3^3 \cdot 11}{5^3}$ which is CM by~$K(\sqrt{-7})$,
\item the point with $\varphi$-coordinate equal to~$-\frac{23549}{5^3}$ which is CM by~$K(\sqrt{-2})$.
\end{enumerate}
\end{prop}

\begin{proof}
The unique pole of~$\varphi$ is known to be CM by~$K(\sqrt{3})$. In order to calculate some
other rational CM-points, we compute the ``modular polynomial''~$\Phi_2(X,Y)$ which is the
algebraic relation between~$\varphi$ and~$\varphi \circ w_2$. It is a symmetric
polynomial of bi-degree~$(9,9)$. The polynomial~$\Phi_2(X,X)$ factorizes into:
\begin{align*} \textstyle
\Phi_2(X,X) =
&\left(X - \textstyle\frac{4752}{125}\right)^2
\left(X + \textstyle\frac{23549}{125}\right)
\left(X^3 - X^2 - 992X - 20736\right) \times \\
&\left(X^3 + \textstyle\frac{95568}{125}X^2 - \textstyle\frac{1212672}{125}X - \textstyle\frac{203493376}{125}\right)^2 \times \\
&\left(X^3 - \textstyle\frac{16752}{125}X^2 - \textstyle\frac{22910208}{15625}X + \textstyle\frac{1126199296}{15625}\right)^2.
\end{align*}
Necessarily, the two rational roots of this polynomial are rational CM points on~$\XX(1)$.

On~$\XX_0(2)$ there is a point~$Q$, defined over~$Q(\sqrt{-7})$, such
that~$\varphi(Q) = \varphi(w_2(Q)) = \frac{4752}{125}$. This means that the
point~$\varphi = \frac{4752}{125}$ is CM by~$K(\sqrt{-7})$.

Above the point~$\varphi = -\frac{23549}{125}$ on~$\XX_0(2)$, there are eight points conjugate to
each other over~$\QQ$, and whose field of definition is a degree eight number field~$M$
such that~$K(\sqrt{-2}) \cdot M$ corresponds to the $2$-ray class field of~$K(\sqrt{-2})$.
Thus, the point~$\varphi = -\frac{23549}{125}$ is necessarily CM by~$K(\sqrt{-2})$.
\end{proof}


\end{document}